\documentclass[12pt,reqno]{amsart}

\usepackage[T1]{fontenc}
\usepackage[margin=1in]{geometry}

\usepackage{amsmath,amssymb,amsthm,mathtools,mathrsfs}
\usepackage[dvipsnames]{xcolor}
\usepackage{graphicx}
\usepackage{tikz}
\usepackage{enumerate}
\usepackage{comment}
\usepackage{mathabx}
\usepackage{cite}
\usepackage{url}
\usepackage[colorlinks, linkcolor=MidnightBlue, citecolor=LimeGreen, urlcolor=blue, hypertexnames=false]{hyperref}
\usepackage{cleveref}
\usepackage[most]{tcolorbox}

\allowdisplaybreaks

\newcommand{\red}[1]{\textcolor{red}{#1}}

\newtheorem{theorem}{Theorem}[section]
\newtheorem{lemma}[theorem]{Lemma}
\newtheorem{proposition}[theorem]{Proposition}
\newtheorem{corollary}[theorem]{Corollary}

\theoremstyle{definition}
\newtheorem{definition}[theorem]{Definition}

\newtheorem{example}[theorem]{Example}

\theoremstyle{remark}
\newtheorem{remark}[theorem]{Remark}

\numberwithin{equation}{section}

\subjclass[2020]{Primary 42C15, 46E22; Secondary 30D15.}

\date{\today}

\keywords{Sampling sequences, Interpolating sequences, Hardy space, Model space, De Branges space, Doubling measure, Beurling–type densities.} 

\title{Sampling and Interpolation for a Non-Doubling Measure under Suitable Conditions}

\author[Anjali and R. K. Srivastava]{Anjali and R. K. Srivastava}

\address{Department of Mathematics, Indian Institute of Technology, Guwahati, India.}

\email{manjali@iitg.ac.in, rksri@iitg.ac.in}

\begin{document}
	
\begin{abstract}
Marzo, Nitzan, and Olsen established sharp density criteria for sampling and interpolation in de Branges spaces with doubling phase measures. We show that these criteria remain valid for a family of Hermite–Biehler functions with non-doubling phase measures when the sequences lie in regions where phase separation is equivalent to Euclidean separation.
\end{abstract}
	
	\maketitle

	\section{Introduction}
Model spaces form a fundamental class of Hilbert spaces of analytic functions  in the study of invariant subspaces of the Hardy space $\mathcal H^2(\mathbb{C}_+)$, where $\mathbb{C}_+ = \{z \in \mathbb{C} : \Im z > 0\}$. By Beurling's theorem, every nontrivial invariant subspace is of the form $\Theta \mathcal H^2(\mathbb{C}_+)$ for some inner function $\Theta$, and the corresponding orthogonal complements
\[\mathcal K_\Theta =\mathcal H^2(\mathbb{C}_+) \ominus \Theta \mathcal H^2(\mathbb{C}_+)\]
are called model spaces. 

De Branges spaces provide an equivalent realization of model spaces in the setting of entire functions which connects Hardy space with  theory of entire functions.

Sampling and interpolation in analytic function spaces have been widely studied, with fundamental results characterized in terms of uniform density conditions \cite{landau1967necessary, seip1992densityFockI, fuhr2007sampling, fuhr2017density}. Along this direction, Marzo, Nitzan, and Olsen \cite{MarzoNitzanOlsen2012} studied these problems in de Branges spaces $H(E)$ under the assumption that the phase $\varphi$ induces a doubling measure $\varphi'(x)\,dx$. They also remarked on the possible necessity of this doubling condition and raised a question of whether such results remain valid without it. They obtain sharp density characterizations which states that sampling sequences correspond to the separated sequences with sufficiently large lower density, while interpolating sequences are separated sequences with sufficiently small upper density.

Motivated by this work, we examine the role of the doubling condition. Since their arguments rely heavily on this assumption, it is a  Marzo, Nitzan, and Olsen \cite{MarzoNitzanOlsen2012} posed a question whether it is merely technical or reflects a deeper structural property. In this paper, we address this question by employing a Hermite-Biehler function with non-doubling phase for which the density conditions for sampling and interpolation still hold, using techniques from multiplier space theory.

The structure of the paper is as follows. In section \ref{sec2} we recall some basic definitions and results. Afterwards, section \ref{sec3} contains the equivalence of sampling and interpolating sequences between a de Branges space and the corresponding model space. We conclude the paper with section \ref{sec4} by presenting an example of Hermite-Biehler function which is associated with a non-doubling measure and the equivalence of the densities and separability notion with respect to two different phase functions.  
 
\section{Preliminaries}\label{sec2}

In this  section, we recall the basic notations and definitions. Throughout the article, we will denote $l^2$ for the space of all square-summable sequences,  and by  $L^2(\mathbb{R})$ we denote  the space of Lebesgue square-integrable functions on $\mathbb{R}$. The space of analytic functions in the upper half-plane $\mathbb{C}_+$, where $z=x+iy$ with $y>0$, will be denoted by  $\mathrm{Hol}(\mathbb{C}_+)$. 
\subsection{Model and de Branges spaces} The Hardy space $\mathcal H^2(\mathbb{C}_+)$ consists of those functions in \( \text{Hol}(\mathbb{C_+}) \) whose boundary values lie in \( L^2(\mathbb{R}) \). Due to these boundary functions $\mathcal H^2(\mathbb{C}_+)$ can be viewed canonically as a closed subspace of \( L^2(\mathbb{R}) \). For detailed study of Hardy spaces one can refer to \cite{Nikolskii2002I,Nikolskii2002II,Mashreghi_representaion_theory}. Every function in the Hardy space admits a canonical decomposition into \emph{inner} and \emph{outer} functions (Section 13.7, \cite{Mashreghi_representaion_theory}). A bounded analytic function in the upper half-plane is called an \emph{inner} function if  \( |u(\xi)| = 1 \) almost everywhere on \( \mathbb{R} \). 
Among the most extensively studied inner functions in
this setting are meromorphic inner functions (MIF) \cite{marakov2005}; that is, if it admits a meromorphic extension to the entire complex plane. It is a well-known result that every MIF can be expressed in the form
\[\Theta(z) =  e^{iaz} \prod_{n \geq 1} \xi_n \frac{z-\lambda_n}{z-\overline{\lambda_n}}, \quad z \in \mathbb{C_+},\]
where  \( a \geq 0 \) is a constant and $\xi_n=\frac{|\lambda^2_n+1|}{\lambda^2_n+1}$, with \(\lambda_n\subset \mathbb{C_+}\) is a sequence which satisfies the Blaschke condition
\[\sum_{n \geq 1} \frac{\Im \lambda_n}{1+|\lambda_n|^2} < \infty,\]see \cite{marakov2005,Mitkovski2010}. Since \( |\Theta(x)|=1 \) for $ x\in\mathbb R $, one can  write $\Theta(x)=e^{i\varphi(x)},$ where $ \varphi:\mathbb R\to\mathbb R $ is a continuous real-valued function, called the argument function or \emph{phase function} of $ \Theta $. We denote the model space corresponding to an inner function $\Theta$ in the upper half-plane by
\[\mathcal K_\Theta = \mathcal H^2(\mathbb{C}_+) \ominus \Theta \mathcal H^2(\mathbb{C}_+),\]which is also a Hilbert space with reproducing kernel 
\begin{equation}\label{reproducing kernel}
	k_\lambda^\Theta(z)=\frac{1}{2\pi i} \frac{1-\overline{\Theta(\lambda)}\Theta(z)}{\overline{\lambda}-z},\qquad \lambda\in \mathbb C^+.
\end{equation}
If $\Theta $ is meromorphic, then all the elements of $\mathcal K_\Theta$ are meromorphic and one can extend equation \eqref{reproducing kernel} to $\lambda\in \mathbb R$.  For detailed study of model spaces, see \cite{mashregi, Nikolskii1986}.

We now study the connection between model spaces and de Branges spaces, which relates Hardy space theory and the theory of entire functions.
A \emph{de Branges space} $H(E)$ is a Hilbert space of entire functions associated with a Hermite--Biehler ($HB$) function $E$, that is, an entire function satisfying 
\[
|E(z)| > |E^\#(z)|, \qquad z \in \mathbb{C}_+,
\]
where $E^\#(z) := \overline{E(\overline{z})}$. The space $H(E)$ is  defined as
\[H(E):=\left\{F-\text{entire}: \frac{F}{E},\frac{F^\#}{E}\in \mathcal H^2(\mathbb{C}_+)\right\},\] with norm $\|F\|_{H(E)}^2 := \int_{\mathbb{R}} \left| \frac{F(x)}{E(x)}\right|^2 \, dx$. In particular, for $E= e^{-i\pi z}$, the space $H(E)$ is called as the \emph{Paley-Wiener space} and we denote this space by $PW$. Every de Branges space is
a reproducing kernel Hilbert space whose reproducing kernel is given by (see \cite{deBranges1968, Levin1996})
\[K_w^E(z) = \frac{\overline{E(w)}\,E(z) - \overline{E^\#(w)}\,E^\#(z)}{2\pi i\,(\overline{w} - z)},\qquad z,w\in \mathbb C.\] 
De Branges spaces and model spaces are related through a unitary mapping $\mathcal U$ from $H(E)$ onto a model space $\mathcal K_\Theta$ which corresponds to a meromorphic inner function $\Theta(z) = \frac{E^\#(z)}{E(z)}$; for $ z\in \mathbb C_+$, (see \cite{Baranov2000}) and is given by  \[
\mathcal  UF = \frac{F}{E},
\]
The phase function $\varphi$ associated with the inner function $\Theta$ is referred to  as the phase function associated with $E$. In this paper, we consider the measure corresponding to phase function $\varphi$ as 
\begin{equation}\label{measure}
\mu_\varphi(I)=\int_I \varphi'(x)dx.
\end{equation}
Additionally, the measure 
\(\mu_\varphi\) is called  a doubling measure if there exists a constant \(C>0\) such that $\mu_\varphi(2I)\le C\,\mu_\varphi(I)$ for every interval \(I\subset\mathbb R\), where \(2I\) denotes the interval with
the same center and twice the length.

\subsection{Sampling, Interpolating sequences and Beurling-type Densities}
We now recall the notions of sampling and interpolating sequences for a general Hilbert space $H$ of functions analytic on a domain $\Omega \subset \mathbb{C}$ with reproducing kernels $k_z$.

\begin{definition}
Let $\Lambda \subset \Omega$ be a sequence of points. Then
\begin{itemize}
\item[(a)] $\Lambda$ is called \emph{sampling} for $H$ if there exist constants $A, B > 0$ such that for all $f \in H$,
\[A\|f\|^2_H \leq \sum_{\lambda \in \Lambda} \frac{|f(\lambda)|^2}{\|k_\lambda\|^2_H} \leq B\|f\|^2_H.\]
\item[(b)] $\Lambda$ is called \emph{interpolating} for $H$ if for every sequence $\{a_\lambda\}$ satisfying\[\sum_{\lambda \in \Lambda} \frac{|a_\lambda|^2}{\|k_\lambda\|^2_H} < \infty,\]
there exists $f \in H$ such that $f(\lambda) = a_\lambda$.
\end{itemize}
\end{definition}
Observe that a sequence $\Lambda$ is sampling if and only if the normalized reproducing kernels $\{\frac{k_\lambda}{\|k_\lambda\|_H}\}$ forms a frame in $H$, while $\Lambda$ is interpolating if and only if it is a Riesz sequence in $H$. For these related notions, see \cite{Seip2004}. 

We use the following notion of separation.
\begin{definition}
A real sequence $\Lambda$ is called $\varphi-$separated  if there exists $\epsilon>0$ such that $|\varphi(\lambda)-\varphi(\lambda')|>\epsilon$ whenever $\lambda,\lambda'\in \Lambda $ and $\lambda\not=\lambda'$.
\end{definition} 
If $\varphi$ is the identity function, then $\Lambda$ is said to be
Euclidean separated. To study the real sampling and interpolating sequences for $H(E)$, we introduce the following Beurling-type densities.
\begin{definition}
Given $\varphi$-separated sequence $\Lambda=\{\lambda_k\}\subset \mathbb R$ and a measure $\mu$ as defined in \eqref{measure}, the lower and upper uniform densities of $\Lambda$ are given by 
\[D_\varphi^-(\Lambda)=\liminf_{r\to\infty} \inf_{I:\mu(I)=r}\frac{|\Lambda \cap I|}{r}, \qquad D_\varphi^+(\Lambda)=\limsup_{r\to\infty} \sup_{I:\mu(I)=r}\frac{|\Lambda \cap I|}{r},\]respectively, where  $I$ are bounded intervals in $\mathbb R$ and $|\Lambda\cap I|$ denotes the number of points in $\Lambda\cap I.$
\end{definition} 

\section{Results for Buliding the Example}\label{sec3}
Throughout this section, we assume  that $U$ and $V$ are meromorphic inner functions and assume that \(\phi\) is an onto multiplier from \(\mathcal K_U\) onto \(\mathcal K_V\); that is, function $\phi$ is analytic in $\mathbb C_+$ which satisfies $
\phi \mathcal K_U=\mathcal K_V.$
Let \(\Gamma=\{\gamma_n\}\subset\mathbb R\), and assume that for $n\ge 1,$ $k^U_{\gamma_n}$ and $ k^V_{\gamma_n}$
are non-zero reproducing kernels of $\mathcal K_U$ and $\mathcal K_V$ respectively. Define the operator $R_U$ from $\mathcal K_U$ onto $R_U(\mathcal K_U)\subset\ell^2$, and similarly define the operator $R_V$ from $\mathcal K_V$ onto $R_V(\mathcal K_V)\subset\ell^2$ as
\begin{equation}\label{operator_R_U}
R_Uf:=\left(\frac{f(\gamma_n)}{\|k^U_{\gamma_n}\|_{\mathcal K_U}}\right)_{n\ge1}\text{ and } R_V  g:=\left(\frac{g(\gamma_n)}{\|k^V_{\gamma_n}\|_{\mathcal K_V}}\right)_{n\ge1}
\end{equation}and the operator $D$ on $\ell^2$ by,
\begin{equation}\label{operator_D}
	D(c_n)=(d_n c_n),
\end{equation}
where $d_n:=\frac{\phi(\gamma_n)\|k^U_{\gamma_n}\|_{\mathcal K_U}}{\|k^V_{\gamma_n}\|_{\mathcal K_V}}$. 
In this section, we establish the relationship between sampling and interpolating sequences in two model spaces, as well as between a model space and its corresponding de Branges space.
\subsection{Sampling and Interpolation in Two Model Spaces}
The following result establishes necessary and sufficient conditions on an onto multiplier $\phi$ between two model spaces under which the sampling property transfers from
one model space to the other.
 
\begin{lemma}\label{sampling_equivalence}
Let $U,V$ be meromorphic inner functions, and assume that \(\phi\) is an onto multiplier from \(\mathcal K_U\) onto \(\mathcal K_V\). Then for any sampling sequence \( \Gamma=\{\gamma_n\}\subset \mathbb R\) for \(\mathcal K_U\), the following are equivalent
\begin{enumerate}[(i)]
\item  $\Gamma \text{ is sampling for }\mathcal K_V;$
\item $D|_{\mathcal R}:\mathcal R\to \ell^2 \text{ is bounded and bounded below},$
\end{enumerate}where the space $ \mathcal R:= R_U( \mathcal K_U)$ with the operator $R_U$ and $D$ are defined as above in equations \eqref{operator_R_U} and \eqref{operator_D} respectively.
\end{lemma}
\begin{proof}
We  observe that for every \(f\in \mathcal K_U\),
\[R_V M_\phi f=\frac{\phi(\gamma_n)f(\gamma_n)}{\|k^V_{\gamma_n}\|_{\mathcal K_V}}=d_n\frac{f(\gamma_n)}{\|k^U_{\gamma_n}\|_{\mathcal K_U}} = D R_U f,\] where $M_\phi$ is a multiplier operator from $\mathcal K_U$ into $\mathcal K_V$ such that $M_\phi f=\phi f$. Hence \begin{equation}\label{R_VM_Phi=DR_U}
	R_V M_\phi=DR_U.
\end{equation}
Since $\phi$ is an onto multiplier,
the operator $M_\phi $ is surjective. Also, \( M_\phi\) is injective. Indeed, if \(\phi f=0\) and \(\phi\not\equiv 0\), then \(f\) becomes  zero on a nonempty open set, hence $f\equiv 0$, since $f$ is analytic. By closed graph theorem, $M_\phi $ is a bounded bijection from $\mathcal K_U$ onto $ \mathcal K_V$. Also, by bounded inverse theorem one can show that its inverse $M_\phi^{-1}$ is bounded.

Now suppose \(\Gamma\) is a sampling set for \(\mathcal K_U\) then \(R_U\) is bounded and bounded below.
Thus \(R_U\) is a bounded bijection from \(\mathcal K_U\) onto the range space of $R_U$ that is, \(\mathcal R\). Also, its inverse \(R_U^{-1}:\mathcal R\to \mathcal K_U\) is bounded, since for $x=R_Uf\in \mathcal R$,
\[\|R_U^{-1}x\|_{\mathcal K_U}=\|f\|_{\mathcal K_U}\le \frac{1}{\sqrt c}\|R_Uf\|_{\ell^2}=\frac{1}{\sqrt c}\|x\|_{\ell^2},\] the inequality here follows from the fact that $\Gamma$ is sampling for $\mathcal K_U$; that is there exists some $c>0 $ such that $c\|f\|_{\mathcal K_U}^2\le\|R_Uf\|_{\ell^2}^2$. Define
\[A:=R_U M_\phi^{-1}:\mathcal K_V\to \mathcal R\]
is a bounded operator, since both $R_U$ and $M_\phi^{-1} $ are bounded operators. This operator is injective also because \(R_U\) and \(M_\phi^{-1}\) are injective, and it is also surjective because if $x\in \mathcal R$, then \(x=R_U f\) for some \(f\in \mathcal K_U\), and with \(g = M_\phi f\in \mathcal K_V\) we have $A g=R_U M_\phi^{-1}(g) = R_U f = x.$ Thus, \(A\) is a bounded bijection. Its inverse is given by $A^{-1}=M_\phi R_U^{-1}:\mathcal R\to \mathcal K_V,$ which is bounded. This implies that $A$ is an isomorphism. From equation \eqref{R_VM_Phi=DR_U} we get, 
\begin{equation}\label{R_V=D A}
	R_V = D|_{\mathcal R} A.
\end{equation}

 We now prove the equivalence. Suppose that \(\Gamma\) is sampling for \(\mathcal K_V\). Then \(R_V\) is bounded above and bounded below. Since \(A^{-1}\) is bounded, \(D|_\mathcal R\) is also bounded.  Also, if \(x\in \mathcal R\), then there exists some $g\in \mathcal K_V$ such that \(x=Ag\). Since \(A\) is an isomorphism, there exists \(C>0\) such that
\[\|Ag\|_{\ell^2}\le C\|g\|_{\mathcal K_V}.\]Here  \(R_V\) is bounded below, therefore $~ \|R_V g\|_{\ell^2}\ge c_0\|g\|_{\mathcal K_V}$ for some $c_0>0$ and for all $ g\in \mathcal K_V$, then \[\|Dx\|_{\ell^2}=\|R_V g\|_{\ell^2}\ge c_0\|g\|_{\mathcal K_V}\ge\frac{c_0}{C}\|x\|_{\ell^2}.\] Therefore, \(D|_\mathcal R\) is bounded below.
 	
Conversely, suppose that \(D|_\mathcal R\) is bounded and bounded below. Since \(A:\mathcal K_V\to \mathcal R\) is an isomorphism and using equation \eqref{R_V=D A} it follows that \(R_V\) is bounded and bounded below. Hence \(\Gamma\) is sampling for \(\mathcal K_V\).
\end{proof}
The following result gives a sufficient condition for a set to be interpolating for one model space when it is interpolating for the other.
\begin{lemma}\label{interpolation_equivalence}
 Let \(U,V\) be meromorphic inner functions, and let $\phi$ be an onto multiplier from $\mathcal K_U $ onto $\mathcal K_V$. Let sequence \(\Gamma=\{\gamma_n\}\subset\mathbb R\). Then $\Gamma$ is interpolating for $\mathcal K_U$ if and only if $\Gamma$ is interpolating for $\mathcal K_V$.
 \end{lemma}
\begin{proof}
 Since $\phi$ is an onto multiplier, so the operator \(M_\phi\) is surjective. Therefore, \(M_\phi\) is invertible and $R_VM_\phi=DR_U$ as proved in equation \eqref{R_VM_Phi=DR_U}. Since \(M_\phi\) is invertible, this implies \(M_\phi^*\); adjoint of the operator $M_\phi$, is also invertible. Also, for $\lambda\in\mathbb R,$ we have
 $M_\phi^*k_\lambda^V=\overline{\phi(\lambda)}k_\lambda^U$ and $M_\phi^*$ is bounded and bounded below. Hence there exist constants \(b,B>0\) such that
 \[b\|k_\lambda^V\|_{\mathcal K_V} \le \|M_\phi^*k_\lambda^V\|_{\mathcal K_U}=|\phi(\lambda)| \|k_\lambda^U\|_{\mathcal K_U}\le B\|k_\lambda^V\|_{\mathcal K_V},\qquad \lambda\in\mathbb R.\]On taking \(\lambda=\gamma_n\),
 \[b\le|\phi(\gamma_n)|\frac{\|k^{U}_{\gamma_n}\|_{\mathcal K_U}}{\|k^{V}_{\gamma_n}\|_{\mathcal K_V}} \le B, \qquad n\ge1,\]
 that is, $b\le |d_n|\le B$ for $ n\ge1,$ where $d_n$ is as defined in \eqref{operator_D}. In particular, $\displaystyle\inf_n |d_n|>0.$
 	
Now assume first that \(\Gamma\) is interpolating for \(\mathcal K_U\) then \(R_U\) is onto \(\ell^2\). Since \(M_\phi\) is onto,
\[R_V(\mathcal K_V)=R_V(M_\phi \mathcal K_U)= DR_U(\mathcal K_U)=D(\ell^2).\]
Since \(\displaystyle\inf_n|d_n|>0\), the diagonal operator \(D\) is onto \(\ell^2\). Hence $R_V(\mathcal K_V)=\ell^2,$ this implies \(\Gamma\) is interpolating for \(\mathcal K_V\).
 	
Conversely, if \(\Gamma\) is interpolating for \(\mathcal K_V\), then by applying same argument  to the inverse multiplier
$M_\phi^{-1}$ one can show that \(\Gamma\) is interpolating for \(\mathcal K_U\).
 \end{proof}
 The following corollary is a consequence of Lemma \ref{sampling_equivalence}, and \ref{interpolation_equivalence}.
 
 \begin{corollary}\label{Invertible_multipliers_preserve_S&I}
Let $U,V$ be meromorphic inner functions, and let $\phi $
 be onto multiplier from $\mathcal K_U$ onto $\mathcal K_V$. Then, for every $\Gamma=\{\gamma_n\}\subset\mathbb R$ the following holds,
 		
\begin{enumerate}[(i)]
 \item $\Gamma$ is interpolating for $\mathcal K_{U}$ if and only if
 			$\Gamma$ is interpolating for $\mathcal K_{V}.$
\item $\Gamma$ is sampling for $\mathcal K_{U}$ if and only if
$\Gamma$ is sampling for $\mathcal K_{V}.$
\end{enumerate}
\end{corollary}
\begin{proof} The proof follows from  Lemma \ref{interpolation_equivalence} and \ref{sampling_equivalence}.
\end{proof}
 	We now study the sampling and interpolation properties under the unitary equivalence between de Branges spaces and model spaces.

\subsection{Sampling and Interpolation in Model and De Branges space}
In this subsection, we show that, under the unitary map between $H(E)$ and $\mathcal K_\Theta$, real sampling and real interpolation are invariant between these two spaces.
		
\begin{proposition}\label{modelspace_de_brange}
Let \(E\in HB\) and  $\Theta$ be the corresponding meromorphic inner function. Then for any sequence $\Gamma=\{\gamma_n\}\subset \mathbb R$ the following holds,
 		
\begin{enumerate}[(i)]
\item $\Gamma$ is sampling for $ H(E)$ if and only if $\Gamma $  is sampling for  $\mathcal K_\Theta$.
\item $\Gamma$  is interpolating for $H(E)$ if and only if	$\Gamma$  is interpolating for $\mathcal K_\Theta.$
\end{enumerate}
\end{proposition}
 	
\begin{proof}Set \[\Theta:=\frac{E^\#}{E},\qquad E^\#(z):=\overline{E(\bar z)}.\]
Then
\[\mathcal U:H(E)\to \mathcal K_\Theta,\qquad \mathcal UF=\frac{F}{E},\]
is a unitary operator \cite{Baranov2000}. Define operator $R_E$ and $R_\Theta $ from $H(E)$ and $\mathcal K_\Theta$ respectively into $\ell^2$ as
\[R_EF:=\left(\frac{F(\gamma_n)}{\|K^E_{\gamma_n}\|_{H(E)}}\right)_{n\ge1},\qquad R_\Theta f:=\left(\frac{f(\gamma_n)}{\|k^\Theta_{\gamma_n}\|_{\mathcal K_\Theta}}\right)_{n\ge1},\] where \(K^E_\lambda\) and \(k^\Theta_\lambda\) are the reproducing kernels of \(H(E)\) and \(\mathcal K_\Theta\), respectively. Assume $D$ is an operator on $\ell^2$ such that
\(D(c_n)=(d_nc_n),\) where $d_n:=\frac{|E(\gamma_n)|}{E(\gamma_n)}$. For \(\lambda\in\mathbb R\), the reproducing kernels are related by
$K^E_\lambda(z)=\overline{E(\lambda)}E(z) k^\Theta_\lambda(z).$ Indeed, for a function \(F\in H(E)\),
\[F(\lambda)=E(\lambda)(\mathcal UF)(\lambda)=E(\lambda)\langle \mathcal UF,k^\Theta_\lambda\rangle_{\mathcal K_\Theta}=\left\langle F,\overline{E(\lambda)}E(\cdot)k^\Theta_\lambda(\cdot)\right\rangle_{H(E)}.\]
Then for each \(n\), \[(R_\Theta \mathcal UF)_n=\frac{F(\gamma_n)/E(\gamma_n)}{\|k^\Theta_{\gamma_n}\|_{\mathcal K_{\Theta}}}=\frac{|E(\gamma_n)|}{E(\gamma_n)}\frac{F(\gamma_n)}{\|K^E_{\gamma_n}\|_{H(E)}}=d_n\,(R_EF)_n,\]
since, $\|K^E_\lambda\|^2_{H(E)}=K^E_\lambda(\lambda)=|E(\lambda)|^2\,k^\Theta_\lambda(\lambda)=|E(\lambda)|^2\|k^\Theta_\lambda\|^2_{\mathcal K(\Theta)}.$ Thus $R_\Theta \mathcal U=DR_E.$ Since \(|d_n|=1\), the diagonal operator \(D\) is unitary.
 		
If \(\Gamma\) is sampling for \(H(E)\), then \(R_E\) is bounded and bounded below. Because \(\mathcal U\) and \(D\) are unitary, the identity $R_\Theta \mathcal U=DR_E$ shows that \(R_\Theta\) is bounded and bounded below. Hence \(\Gamma\) is sampling for \(K_\Theta\). The converse follows by the same argument.
Similarly, since $D$ is onto $\ell^2$ we have,
\[R_\Theta(\mathcal K_\Theta)=R_\Theta(\mathcal UH(E))=DR_E(H(E))=\ell^2.\]
Thus, interpolation is equivalent in $H(E)$ and $\mathcal K_\Theta$, this completes the proof.
\end{proof} 
In the next section, we use a Hermite Behiler function corresponding to a phase function and show that it is not doubling on $\mathbb R $ and we connect the separability and density corresponding to two different phase functions.
 \section{Example }\label{sec4}
 First, we recall few fundamental definitions and results. Let \(\Theta\) be a meromorphic inner function in the upper half-plane \(\mathbb C_+\), with zeros
 $z_n=x_n+i y_n;~ y_n>0,$ which  satisfies the Blaschke condition \[\sum_n \frac{y_n}{1+x_n^2+y_n^2}<\infty.\]
 Up to  unimodular constants, $\Theta$ has the representation
\begin{equation}\label{mif}
 \Theta(z)=e^{iaz}\prod_{n}c_n\frac{z-z_n}{z-\overline{z_n}},\qquad a\geq 0,\quad|c_n|=1.
\end{equation}
The poles of its meromorphic continuation are precisely the reflected points
$\overline{z_n}=x_n-i y_n.$ Since \( |\Theta(x)|=1 \) for \(x\in\mathbb{R}\), we may write $\Theta(x)=e^{i\vartheta(x)},$ where \(\vartheta\) is a continuous increasing phase function. Taking logarithmic derivatives on both sides of equation \eqref{mif} and then differentiate, we get
\[\frac{\Theta'(z)}{\Theta(z)}=ia+\sum_{n}\left(\frac{1}{z-z_n}-\frac{1}{z-\overline{z_n}}\right).\]
Also, for \(x\in\mathbb{R}\),
 \[\frac{\Theta'(x)}{\Theta(x)}=i\vartheta'(x).\]
Therefore, comparing the two equations
above,
\begin{equation}\label{eq:phase_derivative}
	\vartheta'(x)=a+\sum_{n}
	\frac{2y_n}{(x-x_n)^2+y_n^2}.
\end{equation}
 For a detailed study of meromorphic inner functions; see
 \cite{marakov2005, Nikolskii2002I, Nikolskii2002II}.

 \subsection{Non-doubling measure}\label{equivalence} We make use of an example in which $E$ is from Hermite-Behiler class and establish that its corresponding phase
 measure is not doubling on $\mathbb R$. 
\begin{lemma}\label{non_doubling_measure}
 Let \(0<\delta<\tfrac14\), and define
\[E_\delta(z)=(z+i)\prod_{k=1}^{\infty}\left(1-\frac{z}{k-\delta-i k^{-4\delta}}\right)\left(1-\frac{z}{-k+\delta-i k^{-4\delta}}\right).\]
Let \(\varphi_\delta\) be the phase function associated with \(E_\delta\) then the measure $\mu_\delta:=\varphi_\delta'(x)\,dx$ is not doubling on \(\mathbb R\).
\end{lemma}
\begin{proof}
If the Hermite-Behiler function $E_\delta$  has zeros
\(\lambda_n=x_n-i y_n\) where $y_n>0$, in the lower half plane then the zeroes of its corresponding meromorphic inner function \(\Theta_\delta=E_\delta^\#/E_\delta\) in \(\mathbb C_+\) are
$x_n+i y_n$ in the upper half plane and the corresponding phase derivative will be \[\varphi'(x)=\sum_n \frac{2y_n}{(x-x_n)^2+y_n^2},\] up to the addition of the constant term corresponding to an exponential factor; see \eqref{eq:phase_derivative}. Here, there is no such extra constant term, and the zeros of
\(\Theta_\delta=E_\delta^\#/E_\delta\) in \(\mathbb C_+\) are
\[i,\qquad k-\delta+i k^{-4\delta},\qquad -k+\delta+i k^{-4\delta}\qquad (k\ge1).\]
Hence
\begin{equation}\label{eq:phase-derivative}
\varphi_\delta'(x)=\frac{2}{x^2+1}+\sum_{m=1}^\infty\frac{2m^{-4\delta}}{(x-(m-\delta))^2+m^{-8\delta}}+\sum_{m=1}^\infty\frac{2m^{-4\delta}}{(x-(-m+\delta))^2+m^{-8\delta}}.
\end{equation}
We shall show that \(\mu_\delta\) is not a doubling measure.
For each \(k\ge1\), we assume \[x_k:=k-\delta,\qquad I_k:=\left[x_k+\frac18,\;x_k+\frac38\right].\]Then \(I_k\) has center \(x_k+\tfrac14\) and radius \(\tfrac18\), so $2I_k=\left[x_k,\;x_k+\frac12\right].$
We now estimate \(\mu_\delta(2I_k)\) from below and \(\mu_\delta(I_k)\) from above.
 	
From \eqref{eq:phase-derivative}, we only keep \(m=k\) term in the first sum,
we obtain
\[\mu_\delta(2I_k)\ge\int_{2I_k}\frac{2k^{-4\delta}}{(x-x_k)^2+k^{-8\delta}}\,dx.\]
After the change of variable \(t=x-x_k\). Since \(2I_k=[x_k,x_k+\tfrac12]\),
\[\mu_\delta(2I_k)\ge\int_0^{1/2}\frac{2k^{-4\delta}}{t^2+k^{-8\delta}}dt=2\arctan\!\left(\frac{1}{2k^{-4\delta}}\right)\ge2 \arctan\!\left(\frac12\right).\]
Because \(k^{-4\delta}\le1\) for every \(k\ge1\).
Thus there exists a constant \(C_0>0\) which does not depend on \(k\), such that
\begin{equation}\label{eq:lower-bound}
 \mu_\delta(2I_k)\ge C_0 \qquad (k\ge1).
\end{equation}
Therefore, we have shown the lower bound for $\mu_\delta(2I_k)$. 

Now we move towards showing the upper bound for $\mu_\delta(I_k)$. We divide the measure $\mu_\delta(I_k)$ into three parts as 
\[\mu_\delta(I_k)=\int_{I_k}\frac{2dx}{x^2+1}+\sum_{m=1}^{\infty}\int_{I_k}\frac{2m^{-4\delta}}{(x-(m-\delta))^2+m^{-8\delta}}\,dx+\sum_{m=1}^\infty \int_{I_k}\frac{2m^{-4\delta}}{(x-(-m+\delta))^2+m^{-8\delta}}\,dx\]
Suppose $S_1:=\int_{I_k}\frac{dx}{x^2+1}$ and  \[S_2:=\sum_{m=1}^{\infty}\int_{I_k}\frac{m^{-4\delta}}{(x-(m-\delta))^2+m^{-8\delta}}\,dx, \quad S_3:=\sum_{m=1}^\infty \int_{I_k}\frac{m^{-4\delta}}{(x-(-m+\delta))^2+m^{-8\delta}}\,dx \]
First, we will calculate $S_2$. For $m=k$, Again by \eqref{eq:phase-derivative} and change of variables,
\begin{equation}\label{central_bound}
\int_{I_k}\frac{k^{-4\delta}}{(x-x_k)^2+k^{-8\delta}} dx=\int_{1/8}^{3/8}\frac{k^{-4\delta}}{t^2+k^{-8\delta}}\,dt\le k^{-4\delta}\int_{1/8}^{3/8}\frac{dt}{t^2}\le C_1 k^{-4\delta},
\end{equation}where $C_1>0$ and	since \(t\ge \tfrac18\) on the interval of integration. So the central peak contributes at most \(C_1k^{-4\delta}\) on \(I_k\).

For \(m\neq k\). For \(x\in I_k\),
$|x-(m-\delta)|\ge \frac12 |m-k|,$ so \[\int_{I_k}\frac{m^{-4\delta}}{(x-(m-\delta))^2+m^{-8\delta}}\,dx\le|I_k|\frac{m^{-4\delta}}{\frac{1}{4}|m-k|^2}\le c_0\frac{m^{-4\delta}}{|m-k|^2},\qquad c_0>0.\]
Summing over \(m\neq k\),
\begin{equation}\label{positive_bound}
\sum_{m\neq k}\int_{I_k}\frac{m^{-4\delta}}{(x-(m-\delta))^2+m^{-8\delta}}\,dx\le c_0\sum_{m\neq k}\frac{m^{-4\delta}}{|m-k|^2}\le C_2 k^{-4\delta},\qquad C_2>0.
\end{equation}
For the last inequality, we divide the sum into two parts. If $k/2\le m\le 2k$, then $m^{-4\delta}\asymp k^{-4\delta}$, and hence 
\begin{equation}\label{for_m_not_equal_k}
	\sum_{\substack{m\neq k\\ k/2\le m\le 2k}}\frac{m^{-4\delta}}{|m-k|^2}\le c_1 k^{-4\delta}\sum_{j\neq 0}\frac1{j^2} \le c_2 k^{-4\delta}, \qquad c_1,c_2>0.
	\end{equation} 
	If $m<k/2$ or $m>2k$, then $|m-k|\asymp m+k$, so \begin{equation}\label{for_m<k/2}
		\sum_{m<k/2} \frac{m^{-4\delta}}{|m-k|^2}\le c_3\sum_{m<k/2} \frac{m^{-4\delta}}{(m+k)^2}\le\frac{c_3}{k^2}\sum_{m<k/2} m^{-4\delta} \le c_4 k^{-1-4\delta}\le c_4k^{-4\delta},
		\end{equation}where $c_3,c_4>0$ and similarly by integral comparison we get, 
		\begin{equation}\label{for_m>2k}
			\sum_{m>2k} \frac{m^{-4\delta}}{|m-k|^2}\le c_5\sum_{m>2k} \frac{m^{-4\delta}}{(m+k)^2} \le c_5\sum_{m>2k} m^{-2-4\delta}\le c_6 k^{-1-4\delta}\le c_6 k^{-4\delta}.
		\end{equation}
Here $c_5,c_6>0$. Thus combining all equations \eqref{for_m_not_equal_k}, \eqref{for_m<k/2} and \eqref{for_m>2k}  we get the last inequality of equation \eqref{positive_bound}. Therefore, the total contribution of all positive peaks with \(m\neq k\) is bounded by \(C_2 k^{-4\delta}\), where the constant $C_2 $ depends on $\delta$ only. Thus, using \eqref{central_bound} and \eqref{positive_bound} we get, 
\begin{equation}\label{for_positive}
	S_2=\sum_{m=1}^{\infty}\int_{I_k}\frac{m^{-4\delta}}{(x-(m-\delta))^2+m^{-8\delta}}\,dx\le C k^{-4\delta},\qquad C>0.
\end{equation}

Now, we calculate $S_3$. For \(x\in I_k\) and \(m\ge1\), $|x-(-m+\delta)|\ge k+m-\frac12,$ so
\[\int_{I_k} \frac{m^{-4\delta}}{(x-(-m+\delta))^2+m^{-8\delta}}\,dx\le C_3\frac{m^{-4\delta}}{(k+m)^2}. \]
Hence
\begin{equation}\label{negative_bound}
S_3=\sum_{m=1}^\infty \int_{I_k}\frac{m^{-4\delta}}{(x-(-m+\delta))^2+m^{-8\delta}}\,dx\le C_3\sum_{m=1}^\infty \frac{m^{-4\delta}}{(k+m)^2}\le C' k^{-4\delta},
\end{equation}
where $C_3, C'>0$ and the last inequality can be solved as in a manner similar to the arguments used in  equation \eqref{for_m<k/2} and \eqref{for_m>2k} by dividing the sum into two parts when $m\le k$ and $m> k$. Hence again $S_3$ has \(O(k^{-4\delta})\) because \(4\delta<1\).
Also for $S_1$, when  \(x\in I_k\), we have \(x\asymp k\), so
\begin{equation}\label{bound_I_k}
S_1=\int_{I_k}\frac{dx}{x^2+1}\le \frac{C'''}{k^2},\qquad C'''>0.
\end{equation}
Since \(2>4\delta\), this is also \(O(k^{-4\delta})\).
Combining \eqref{for_positive}, \eqref{negative_bound} and \eqref{bound_I_k}, we obtain
$\mu_\delta(I_k)\le C_{\delta} k^{-4\delta},$  where $C_{\delta}>0$ depends only on $\delta$. Together with \eqref{eq:lower-bound}, this gives
\[\frac{\mu_\delta(2I_k)}{\mu_\delta(I_k)}\ge\frac{C_0}{C_{\delta} k^{-4\delta}}=\frac{C_0}{C_{\delta}}\,k^{4\delta}\xrightarrow{}\infty,\]as $k \to \infty$.	Therefore no uniform doubling constant can exist for \(\mu_\delta\). Hence the measure \(\varphi_\delta'(x)\,dx\) is not doubling on \(\mathbb R\).
\end{proof}
For $z \in \mathbb C_+$ the inner function $	\Theta_\pi(z)=e^{2\pi i z}$ has the corresponding phase function $\varphi_\pi$ as the function $2\pi x$. Therefore, the measure $\mu_\pi$ corresponding to the phase function  is a constant multiple of  Lebesgue measure. Consequently, $\mu_\pi$ is a doubling measure on $\mathbb R$. 

Let $\Theta_\delta=E_\delta^\#/E_\delta$ and $\Psi_\delta(z)=e^{i\pi z}E_\delta(z)$; where $E_\delta$ is as defined in Lemma \ref{non_doubling_measure}. By the standard correspondence between model spaces and de Branges spaces, multiplication by $E_\delta$ maps $\mathcal K_{\Theta_\delta}$ onto $ H(E_\delta)$; see \cite{Baranov2000}. Moreover, it was shown in \cite{LyubarskiiSeip2002} that $E_\delta$ is a Hermite–Biehler function and that $ H(E_\delta)=PW$ with equivalent norms, where $PW$ is the Paley-Wiener space defined in Section \ref{sec2}. Therefore, we obtain $E_\delta\mathcal K_{\Theta_\delta}=PW$. Multiplying this identity by $e^{i\pi z}$ gives $\Psi_\delta\mathcal K_{\Theta_\delta}=e^{i\pi z}PW$. By the standard identification $e^{i\pi z}PW=\mathcal K_{\Theta_\pi}$, it follows that $\Psi_\delta\mathcal K_{\Theta_\delta}=\mathcal K_{\Theta_\pi}$. Hence multiplication by $\Psi_\delta$ maps $\mathcal K_{\Theta_\delta}$ onto $\mathcal K_{\Theta_\pi}$, and therefore $\Psi_\delta$ is an onto multiplier.
 
Although \(\Psi_\delta\) is unbounded in \(\mathbb C_+\), this does not affect the invertibility of the multiplication operator
\[M_{\Psi_\delta}:\mathcal K_{\Theta_\delta}\to \mathcal K_{\Theta_\pi}.\]
Indeed, invertibility depends on boundedness of \(M_{\Psi_\delta}\) as an operator between the Hilbert spaces \(\mathcal K_{\Theta_\delta}\) and \(\mathcal K_{\Theta_\pi}\), not on boundedness of the symbol \(\Psi_\delta\) itself. In the present situation,
\[\|\Psi_\delta f\|_{\mathcal K_{\Theta_\pi}}\asymp \|f\|_{\mathcal K_{\Theta_\delta}},
 \qquad f\in\mathcal K_{\Theta_\delta},\]
so \(M_{\Psi_\delta}\) is bounded and bounded below. Since also
 \[\Psi_\delta\mathcal K_{\Theta_\delta}=\mathcal K_{\Theta_\pi},\]
 the operator \(M_{\Psi_\delta}\) is a bounded bijection, hence invertible.

Since $\Psi_\delta$ is an onto multiplier from $\mathcal K_{\Theta_\delta}$ onto $\mathcal K_{\Theta_\pi}$, Corollary~\ref{Invertible_multipliers_preserve_S&I} implies that
\[\Gamma \text{ is sampling/interpolating for } \mathcal K_{\Theta_\delta}\iff\Gamma \text{ is sampling/interpolating for } \mathcal K_{\Theta_\pi}.\] Moreover, from Proposition \ref{modelspace_de_brange} implies,
\begin{align*}
\Gamma \text{ is sampling/interpolating for } \mathcal K_{\Theta_\delta}
&\iff\Gamma \text{ is sampling/interpolating for } H(E_\delta),\\
\Gamma \text{ is sampling/interpolating for } \mathcal K_{\Theta_\pi}
&\iff\Gamma \text{ is sampling/interpolating for } H(e^{-i\pi z}).
\end{align*}
Hence,
\begin{equation}\label{equivalences}
	\Gamma \text{ is sampling/interpolating for } H(E_\delta)\iff
\Gamma \text{ is sampling/interpolating for } H(e^{-i\pi z}).
\end{equation}
Up to this point, we have established only the equivalence of sampling and interpolating sequences. To verify that $E_\delta$ satisfies the hypothesis of the main result in \cite{MarzoNitzanOlsen2012}, it remains to investigate the relationship between the density and separation properties associated with doubling and non-doubling phase functions.

\subsection{Density and Separability Comparison}
In this subsection, we compare the notions of density and separation induced by the two phase functions $\varphi_\delta$ and $\varphi_\pi$. Here, we also determine the regions in which these notions are equivalent and examine how the asymptotic relation between the phases affects the corresponding upper and lower densities.
\begin{lemma}\label{boundforF_k}
	For $x>0$ and $0<\delta<\frac{1}{4}$, there exists a constant $C>0$ such that
	\[\sum_{\substack{k\ge 1\\ |x-a_k|\ge1}}\left(\frac{y_k}{|x-a_k|}+\frac{y_k}{x+a_k}\right)\leq C\qquad \text{ as } x\to \infty,\]where $y_k=k^{-4\delta}$ and $a_k=k-\delta$. 
	\end{lemma}
	\begin{proof}
		First, assume \[S=\sum_{\substack{k\ge 1\\ |x-a_k|\ge1}}\frac{y_k}{|x-a_k|} \] and divide this sum into three parts when $k\leq\frac{x}{2}$, when $\frac{x}{2}<k\leq 2x$ and when $k>2x$. 
		
		When $k\leq\frac{x}{2}$ then using  $|x-a_k|=|x-k+\delta|=x-k+\delta\ge x-\frac{x}{2}+\delta\ge \frac{x}{2}$ we obtain,
		\begin{align}
			\sum_{\substack{|x-a_k|\ge1\\1\leq k\leq \frac{x}{2}}}\frac{y_k}{|x-a_k|}&=\sum_{\substack{|x-a_k|\ge1\\1\leq k\leq \frac{x}{2}}}\frac{k^{-4\delta}}{|x-a_k|}\leq \sum_{\substack{|x-a_k|\ge1\\1\leq k\leq \frac{x}{2}}}\frac{2k^{-4\delta}}{x}=\frac{2}{x}\sum_{\substack{|x-a_k|\ge1\\1\leq k\leq \frac{x}{2}}} k^{-4\delta}\nonumber\\
			&< \frac{2}{x(1-4\delta)}\left(\frac{x}{2}\right)^{1-4\delta}= C_1 x^{-4\delta}\label{k<=x/2},
		\end{align}
		where $C_1>0$ and the last inequality \eqref{k<=x/2} follows using the fact that $\sum_{k=1}^{M}k^{-p}<\frac{M^{1-p}}{1-p}$, for $0<p<1$. 
		
		When $\frac{x}{2}<k\leq 2x$, then for $0<p<1$ we have $k^{-p}\leq 2^{p} x^{-p}$ and also the following inequality holds 
		\begin{equation}\label{boundsofx-a_k}
			-x+\delta\leq x-a_k \leq \frac{x}{2}+\delta\leq 2x+\delta\leq \lceil 2x+\delta\rceil,
		\end{equation} 
		where $ \lceil \cdot \rceil $ denotes the ceiling function.  Notice that all $a_k's$ are at a distance of 1 from each other therefore between any region such that $j\leq |x-a_k|\leq j+1$ there can be atmost two $a_k's$. Thus one can write, 
		\begin{equation}\label{bound2/j}
			\sum_{\substack{k\ge 1\\j\leq |x-a_k|\leq j+1}}\frac{1}{|x-a_k|}\leq \frac{2}{j}.
		\end{equation} 
		Using above equations \eqref{boundsofx-a_k} and \eqref{bound2/j} we obtain,
		\begin{align}
			\sum_{\substack{|x-a_k|\ge1\\ \frac{x}{2}\leq k\leq 2x}}\frac{k^{-4\delta}}{|x-a_k|}&\leq 2^{4\delta} x^{-4\delta} \sum_{1\leq |x-a_k|\leq\lceil 2x+\delta\rceil}\frac{1}{|x-a_k|}\nonumber\\
			&\leq 2^{4\delta} x^{-4\delta} 2\left(\frac11+\frac12+\frac13+\cdots+\frac{1}{\lceil 2x+\delta\rceil}\right)\nonumber\\
			&\leq 2^{4\delta} x^{-4\delta}2(1+\log (\lceil 2x+\delta\rceil)) \leq 2^{4\delta}10  x^{-4\delta}\log(2+x)\label{loginequality}\\
			&= C_2 x^{-4\delta}\log(2+x), \label{kinbetween}
		\end{align}
		where $C_2>0$ and the last inequality \eqref{loginequality} follows due to $\sum_{j=1}^M\frac{1}{j}\leq 1+\log M. $
		
		When $k> 2x$, then it is easy to verify that $\frac{k}{4}<\frac{k}{2}-\delta<k-\delta-x=|x-a_k|$ and  take
		\begin{align}
			\sum_{\substack{|x-a_k|\ge1 \nonumber\\ k>2x}}\frac{k^{-4\delta}}{|x-a_k|}&\leq \sum_{\substack{|x-a_k|\ge1 \nonumber\\ k>2x}}\frac{4k^{-4\delta}}{k}< 4 \sum_{k>2x} k^{-1-4\delta} \nonumber\\
			&\leq c_1 \int_{2x}^{\infty} t^{-1-4\delta}dt= c_1 \frac{(2x)^{-4\delta}}{4\delta}\nonumber\\
			&\leq C_3 x^{-4\delta} \label{whenk>2x},
		\end{align}
	where $c_1,C_3>0$.	Combining equations \eqref{k<=x/2}, \eqref{kinbetween} and \eqref{whenk>2x},
		\begin{align}
			S&=\sum_{\substack{k\ge 1\\ |x-a_k|\ge1}}\frac{y_k}{|x-a_k|}= \sum_{\substack{|x-a_k|\ge1\\1\leq k\leq \frac{x}{2}}}\frac{y_k}{|x-a_k|}+\sum_{\substack{|x-a_k|\ge1\\ \frac{x}{2}\leq k\leq 2x}}\frac{k^{-4\delta}}{|x-a_k|}+\sum_{\substack{|x-a_k|\ge1\\ k>2x}}\frac{k^{-4\delta}}{|x-a_k|}\nonumber\\
			&\leq C_1 x^{-4\delta}+C_2 x^{-4\delta}\log(2+x)+C_3 x^{-4\delta} \to 0 \qquad \text{ as } x\to \infty. \label{S}
		\end{align}
	Now consider, \begin{align}\sum_{\substack{k\ge 1\\ |x-a_k|\ge1}}\left(\frac{y_k}{|x-a_k|}+\frac{y_k}{x+a_k}\right)&=\sum_{\substack{k\ge 1\\ |x-a_k|\ge1}}\frac{y_k}{|x-a_k|}+\sum_{\substack{k\ge 1\\ |x-a_k|\ge1}}\frac{k^{-4\delta}}{x+k-\delta}\nonumber\\
			&\leq S+ \sum_{k\ge1}  \frac{k^{-4\delta}}{x+k-\delta}\leq S+ \sum_{k\ge1}  \frac{k^{-4\delta}}{k-\delta}\nonumber\\
			& \leq S+ C_4 \label{limitcomparisonbound},
		\end{align}
	where $C_4>0$ and the inequality \eqref{limitcomparisonbound} follows because $\sum_{k\ge1}  \frac{k^{-4\delta}}{(k-\delta)}$ is a convergent series by limit comparison test, on comparing with convergent series $\sum_{k\ge 1}\frac{1}{k^{1+4\delta}}$. Hence, using equation \eqref{S} in \eqref{limitcomparisonbound} we get our required result.
	\end{proof}
	
The next theorem shows that the absolute value of the difference of two phase functions is bounded by a constant term at infinity. Consequently, this bounded perturbation does not affect the upper or lower phase densities of any discrete sequence. In the following theorems we use the notation $f(x)=O(g(x))$ to mean that there exists a constant $C>0$ and $x_0$ such that $|f(x)|\leq C|g(x)|$ for all $x\ge x_0.$

\begin{theorem}\label{densities}
Let $\varphi_{\pi} $ and $\varphi_{\delta}$ be the phase function corresponding to $e^{-i\pi z}$ and $E_\delta$ respectively; for definition of $E_\delta$ see Lemma \ref{non_doubling_measure}. Then
\[\varphi_\delta(x)=\varphi_\pi(x)+O(1)=2\pi x+O(1),\qquad |x|\to\infty.\]
Consequently, for every discrete sequence \(\Gamma\subset\mathbb R\),
\[D_{\varphi_\delta}^{+}(\Gamma)=D_{\varphi_\pi}^{+}(\Gamma),\qquad D_{\varphi_\delta}^{-}(\Gamma)=D_{\varphi_\pi}^{-}(\Gamma).\]
\end{theorem}
\begin{proof}
For \(\Theta_\pi(z)=e^{2\pi i z}\), we have
\[\varphi_\pi(x)=2\pi x.\]
The zeros of \(I_\delta\) in \(\mathbb C_+\) are
\[i,\qquad k-\delta+i k^{-4\delta},\qquad -k+\delta+i k^{-4\delta},\qquad k\geq 1.\]
Hence the phase derivative of \(I_\delta\) is given by, as calculated in Lemma \ref{non_doubling_measure},
\[\varphi_\delta'(x)=\frac{2}{x^2+1}+\sum_{k\geq 1}\left[\frac {2y_k}{(x- a_k)^2+y_k^2}+\frac{2y_k}{(x+a_k)^2+y_k^2}\right],\]
where $a_k= k-\delta$ and $y_k= k^{-4\delta}$.
  We will first consider $x>0$. On integrating from $0$ to $x$, we get 
\[ \varphi_\delta(x)-\varphi_\delta(0)=2\arctan x + 2\sum_{k\ge 1} F_k(x),\]
where \[
	F_k(x):=\arctan\left(\frac{x-a_k}{y_k}\right)+\arctan \left(\frac{x+a_k}{y_k}\right).\]
Now, we consider $F_k(x)$ in two cases, when \textbf{$ a_k\leq x$}, using $\arctan t+ \arctan \frac{1}{t}=\frac{\pi}{2}$ for $t>0$, one can write \[F_k(x)= \pi - \arctan \left(\frac{y_k}{x-a_k}\right)-\arctan \left(\frac{y_k}{x+a_k}\right).\]
Here $a_k\leq x$, $\textbf{1}_{\{k-\delta\leq x\}}=1$. Therefore when $a_k\leq x$, 
\begin{equation}\label{case_a_k<x}
	|F_k(x)-\pi\cdot\textbf{1}_{\{k-\delta\leq x\}} |= \arctan \left(\frac{y_k}{x-a_k}\right)+\arctan \left(\frac{y_k}{x+a_k}\right) \leq \frac{y_k}{x-a_k}+\frac{y_k}{x+a_k}.
\end{equation}
Also, Using the middle expression above we have,
\begin{equation}\label{pi_bound 1}
	|F_k(x)-\pi|\leq \pi \qquad \text{ when } a_k\leq x.
	\end{equation}
When $0< x<a_k$, here again using the $\arctan$ property as in the previous case, we obtain \[F_k(x)=\arctan \left(\frac{y_k}{a_k-x}\right)-\arctan \left(\frac{y_k}{x+a_k}\right).\]
Since here $x<a_k$, $\textbf{1}_{\{k-\delta\leq x\}}=0$. We get,
\begin{equation}\label{case_a_k>x}
	|F_k(x)-\pi\cdot\textbf{1}_{\{k-\delta\leq x\}} |=|F_k(x)-\pi(0)|= |F_k(x)| \leq \frac{y_k}{a_k-x}+\frac{y_k}{x+a_k},
\end{equation}
here also using expression for $F_k(x)$ in the middle, we get 
\begin{equation}\label{pi_bound 2}
	|F_k(x)-\pi(0)|= |F_k(x)| \leq \pi \qquad \text{ when } x<a_k.
	\end{equation}On combining equations \eqref{case_a_k<x} and \eqref{case_a_k>x} we get, for $x>0$
\begin{equation}\label{F_k(x)-pi}
	|F_k(x)-\pi\textbf{1}_{\{k-\delta\leq x\}}|\leq \frac{y_k}{|x-a_k|}+\frac{y_k}{x+a_k}, 
	\end{equation}
and also merging equation \eqref{pi_bound 1} and \eqref{pi_bound 2} we get, 
\begin{equation}\label{complete pi bound}
	|F_k(x)-\pi\textbf{1}_{\{k-\delta\leq x\}}|\leq \pi.
	\end{equation}
	Consider $N(x):=\#\{k\in \mathbb N: k-\delta \leq x\}=\lfloor x+ \delta \rfloor$, where $\lfloor \cdot \rfloor$ denotes Greatest integer function. Therefore, \[ |N(x)-x|\leq 2.\]Additionally, one can rewrite this function as $\displaystyle N(x)=\sum_{k\ge1 }\textbf{1}_{\{k-\delta\leq x\}}.$ Next compute, 
	\begin{align}
	\sum_{k\ge 1}|F_k(x)-\pi\textbf{1}_{\{k-\delta\leq x\}}|&= \sum_{\substack{k\ge 1\\ |x-a_k|<1}}|F_k(x)-\pi\textbf{1}_{\{k-\delta\leq x\}}|+\sum_{\substack{k\ge 1\\ |x-a_k|\ge1}}|F_k(x)-\pi\textbf{1}_{\{k-\delta\leq x\}}|\nonumber\\
	&\leq 2\pi +\sum_{\substack{k\ge 1\\ |x-a_k|\ge1}}\left(\frac{y_k}{|x-a_k|}+\frac{y_k}{x+a_k}\right)\label{sum{F_k-piNx}}.
\end{align}
The last inequality \eqref{sum{F_k-piNx}} follows using equations \eqref{F_k(x)-pi}, \eqref{complete pi bound} and in light of the fact that $a_{k+1}-a_k=1$, therefore there are atmost two centers satisfying $|x-a_k|<1$. Moreover, due to Lemma \ref{boundforF_k} there exists a constant $C>0 $ such that 
\begin{equation}\label{xtoinfinity}
	\sum_{k\ge 1}|F_k(x)-\pi\textbf{1}_{\{k-\delta\leq x\}}|\leq 2\pi +C \qquad \text{ as } x\to \infty.
\end{equation}
Now using equation \eqref{xtoinfinity} proceed to calculate, 
\begin{align}
	|\varphi_\delta(x)-\varphi_\pi(x)|&\leq|2\pi N(x)-2x\pi|+|\varphi_\delta(x)-2\pi N(x)|\leq 4\pi+ |\varphi_\delta(x)-2\pi N(x)|\nonumber\\
	&\leq  4\pi +|2\arctan x+ \varphi_\delta (0)|+ |2\sum_{k\ge 1} F_k(x)-2\pi N(x)|\nonumber\\
	&\leq 5\pi + \varphi_\delta(0) +2|\sum_{k\ge 1} F_k(x)-\pi \textbf{1}_{\{k-\delta\leq x\}}|\nonumber\\
	&\leq C' +2\sum_{k\ge 1}|F_k(x)-\pi\textbf{1}_{\{k-\delta\leq x\}}|\leq C'' +4\pi
\qquad \text{ as } x\to \infty, \label{positive}
\end{align}
where $C''=C+C'=C+5\pi + \varphi_\delta(0)>0$.
  
On the other hand, when $x<0$, observe that $\varphi_\delta'(x)$ is an even function. For $t>0$ we have, 
\[\varphi_\delta(t)-\varphi_\delta(0)=\int_0^t\varphi_\delta'(u)du.\] Similarly we have, \[\varphi_\delta(-t)-\varphi_\delta(0)=\int_0^{-t}\varphi_\delta'(u)du=-\int_0^{t}\varphi_\delta'(-v)dv=-(\varphi_\delta(t)-\varphi_\delta(0)) .\] This implies, 
\[\varphi_\delta(-t)=2\varphi_\delta(0)-\varphi_\delta(t).\]
Now take
\begin{align}
	|\varphi_\delta(-t)+2\pi t|&= |2\varphi_\delta(0)-\varphi_\delta(t)+2\pi t|\leq|2\varphi_\delta(0)| + |\varphi_\delta(t)-2\pi t|\nonumber\\
	&\leq |2\varphi_\delta(0)| + C''+2\pi \qquad \text{ as } t\to \infty. \label{negative}
	\end{align}
Merging both equations \eqref{positive} and \eqref{negative} one can write for large $x\in \mathbb  R$ there exists a constant \(C_0>0\) such that
\[|\varphi_\delta(x)-\varphi_\pi(x)|\leq C_0.\]
	
Let \(J\subset\mathbb R\) be an interval of \(\varphi_\delta\)-length \(r\) i.e. \[J:=\{x\in \mathbb R: a\leq \varphi_\delta(x)\leq a+r\},\] since $\varphi_\delta$ is an increasing function on $\mathbb R$, $J$ will also be an interval on the real line. Take any $x\in J$, this will imply $a\leq \varphi_\delta(x)\leq a+r $. Additionaly,  for large $x$ we have already obtained,
\[|\varphi_\delta(x)-\varphi_\pi(x)|\leq C.\]
This gives \[a-C\leq \varphi_\delta(x)-C\leq \varphi_\pi(x)\leq C+\varphi_\delta(x)\leq C+a+r,\]
which means $\varphi_\pi(x)\in [a-C,a+r+C]$. Therefore, 
\[J\subset \{x: \varphi_\pi(x)\in [a-C,a+r+C]\}.\]
Thus, every $\varphi_\delta$-interval of length $r$ is contained in a $\varphi_\pi$-interval of length $r+2C$. Likewise, one can show
\[\{x: \varphi_\pi(x)\in [a+C,a+r-C]\}\subset J.\] Consequently, \[\{x:\varphi_\pi(x)\in [a+C,a+r-C]\}\subset J\subset\{x:\varphi_\pi(x)\in [a-C,a+r+C]\}.\]
So for every discrete sequence \(\Gamma\subset\mathbb R\),
\[|\bigl(\Gamma\cap \varphi_\pi^{-1}([a+C,a+r-C])\bigr)|\leq|(\Gamma\cap J)|
\leq|\bigl(\Gamma\cap \varphi_\pi^{-1}([a-C,a+r+C])\bigr)|.\]
Dividing by $r$ we obtain,
\begin{equation}\label{densityinequality}
	\frac{(r-2C)|\Gamma\cap \varphi_\pi^{-1}([a+C,a+r-C])|}{r(r-2C)}\leq \frac{|(\Gamma\cap J)|}{r}\leq \frac{(r+2C)|\Gamma\cap \varphi_\pi^{-1}([a-C,a+r+C])|}{r(r+2C)}. 
\end{equation}
Taking supremum over $a\in \mathbb R$ gives
\begin{align*}\left(1-\frac{2C}{r}\right)\sup_{a\in \mathbb R}\frac{|\Gamma\cap \varphi_\pi^{-1}([a+C,a+r-C])|}{(r-2C)}&\leq \sup_{a\in\mathbb  R}\frac{|(\Gamma\cap J)|}{r}\\
	&\leq \left(1+\frac{2C}{r}\right)\sup_{a\in\mathbb  R}\frac{|\Gamma\cap \varphi_\pi^{-1}([a-C,a+r+C])|}{(r+2C)}.
\end{align*}
After taking $\limsup_{r\to\infty}$, we conclude 
\[D_{\varphi_\pi}^{+}(\Gamma)\leq D_{\varphi_\delta}^{+}(\Gamma)\leq D_{\varphi_\pi}^{+}(\Gamma).\]
 Hence
\[D_{\varphi_\delta}^{+}(\Gamma)=D_{\varphi_\pi}^{+}(\Gamma).\]
Similarly, first taking infimum over $a\in \mathbb R$ in \eqref{densityinequality} and then taking $\liminf_{r\to\infty}$ we obtain,
\[D_{\varphi_\delta}^{-}(\Gamma)=D_{\varphi_\pi}^{-}(\Gamma).\]
This proves the result.
\end{proof}
The equality of the phase densities does not imply equivalence of the  corresponding separation conditions. In the next example we see that  the two separation conditions corresponding to $\varphi_\delta$  and $\varphi_\pi$ are not equivalent while the corresponding Beurling densities may agree on a large scale.
\begin{example}
Let \[\Theta_\pi(z)=e^{2\pi i z},\qquad z\in\mathbb C_+,\]
and let \(\Theta_\delta=E_\delta^\#/E_\delta\), where \(E_\delta\) is defined in Lemma \ref{non_doubling_measure}. Let \(\varphi_\pi\) and \(\varphi_\delta\) be the corresponding phase functions. Then \(\varphi_\pi\)-separation and \(\varphi_\delta\)-separation are not comparable.
\end{example}
\begin{proof}
Let \(a_k=k-\delta\) and \(y_k=k^{-4\delta}\). Consider
\[\Gamma_1=\left\{k+\frac13,\ k+\frac23:\ k\geq 2\right\}.\]
Then \(\Gamma_1\) is \(\varphi_\pi\)-separated, because
\(\varphi_\pi(x)=2\pi x\) and \[\left|\varphi_\pi(k+2/3)-\varphi_\pi(k+1/3)\right|=2\left|\frac{2\pi}{3}-\frac{\pi}{3}\right|=\frac{2\pi}{3}>0.\]

However, on the intervals \([k+1/3,k+2/3]\), since $\varphi_\delta'(x)>0$, the phase function $\varphi_\delta$ is strictly increasing, hence one can write
\begin{equation}\label{difference phi}
	\left|\varphi_\delta\left(k+\frac23\right)-\varphi_\delta\left(k+\frac13\right)\right|=\varphi_\delta\left(k+\frac23\right)-\varphi_\delta\left(k+\frac13\right)=\int_{k+\frac13}^{k+\frac23} \varphi_\delta'(x)dx,
\end{equation}
here \[\varphi_\delta'(x)=\frac{2}{x^2+1}+\sum_{j\geq 1}\left[\frac {2j^{-4\delta}}{(x- a_j)^2+j^{-8\delta}}+\frac{2j^{-4\delta}}{(x+a_j)^2+j^{-8\delta}}\right],\] where $a_j=j-\delta$. It is easy to verify that for $k+\frac13\leq x\leq k+\frac23,$ and for all $j\ge 1$, $|x-a_j|\ge \frac{1}{24}(1+|j-k|)$. This gives, 
\begin{equation}
	\sum_{j=1}^{\infty} \frac{j^{-4\delta}}{(x-a_j)^2+j^{-8\delta}}
	\leq \sum_{j=1}^{\infty} \frac{j^{-4\delta}}{(x-a_j)^2}
	\leq 24 \sum_{j=1}^{\infty} \frac{j^{-4\delta}}{(1+|j-k|)^2}\label{24}.
\end{equation}
Note that $\frac{k}{j}\leq 1$ for $j\ge k$  and  we can write $\frac{k}{j}=1+\frac{k-j}{j}\leq 1+ (k-j)$ for $j\leq k$. Therefore, on merging both we deduce 
\[\frac{k}{j}\leq (1+|j-k|).\] 
Using above inequality in equation \eqref{24} we get 
\begin{align}
	\sum_{j=1}^{\infty} \frac{j^{-4\delta}}{(x-a_j)^2+j^{-8\delta}}
		&\leq 24 \sum_{j=1}^{\infty} \frac{k^{-4\delta}(1+|j-k|)^{4\delta}}{(1+|j-k|)^2}= 24 k^{-4\delta}\sum_{j=1}^{\infty} \frac{1}{(1+|j-k|)^{2-4\delta}}\nonumber\\
	&\leq 24 k^{-4\delta} \left(1+ 2 \sum_{m=1}^{\infty}\frac{1}{(1+m)^{2-4\delta}}\right)\nonumber\\
	&\leq C' k^{-4\delta} \label{negativeside},
\end{align}
where $C'>0$ and the last inequality \eqref{negativeside} holds because $\sum_{m=1}^{\infty}\frac{1}{(1+m)^{2-4\delta}}$ is a convergent series. 

Next, considering the positive term series, notice for $k+\frac13\leq x\leq k+\frac23,$ we have $\frac{3(j+k)}{4}\leq x+a_j$, therefore
\begin{align}
	\sum_{j=1}^{\infty} \frac{j^{-4\delta}}{(x+a_j)^2+j^{-8\delta}}
	&\leq \sum_{j=1}^{\infty} \frac{j^{-4\delta}}{(x+a_j)^2}\leq \frac{16}{9}\sum_{j=1}^{\infty} \frac{j^{-4\delta}}{(k+j)^2}= \frac{16}{9} \left(\sum_{j=1}^{k} \frac{j^{-4\delta}}{(k+j)^2}+\sum_{j=k+1}^{\infty} \frac{j^{-4\delta}}{(k+j)^2}\right)\nonumber\\
	&\leq \frac{16}{9} \left(\sum_{j=1}^{k} \frac{j^{-4\delta}}{k^2}+\sum_{j>k}^{\infty} \frac{j^{-4\delta}}{j^2} \right)\leq \frac{16}{9} \left( \frac{1}{k^2}\sum_{j=1}^{k}j^{-4\delta}+ \int_{k}^{\infty}t^{-4\delta-2}dt\right)\nonumber\\
	&\leq \frac{16}{9} \left(\frac{k^{1-4\delta}}{k^2(1-4\delta)}+ \frac{k^{-1-4\delta}}{1+4\delta}\right)= C'' k^{-4\delta-1}, \label{positiveside}
\end{align}
where $C''>0$ and we use   $\sum_{k=1}^{M}k^{-p}<\frac{M^{1-p}}{1-p}$; where $0<p<1$, for the last inequality \eqref{positiveside}.

Also, for  $k+\frac13\leq x\leq k+\frac23,$ we have 
\begin{equation}\label{boundfor x}
	\frac{1}{1+x^2}\leq \frac{1}{x^2}\leq \frac{1}{k^2}.
\end{equation}
Thus, using equations \eqref{negativeside}, \eqref{positiveside}, \eqref{boundfor x} for equation \eqref{difference phi} we get, 
\begin{align*}
	\left|\varphi_\delta\left(k+\frac23\right)-\varphi_\delta\left(k+\frac13\right)\right|&=\int_{k+\frac13}^{k+\frac23}\left(\frac{1}{x^2+1}+\sum_{j\geq 1}\left[\frac {j^{-4\delta}}{(x- a_j)^2+j^{-8\delta}}+\frac{j^{-4\delta}}{(x+a_j)^2+j^{-8\delta}}\right]\right)dx\\
	&\leq C' k^{-4\delta} + C''k^{-4\delta-1} + k^{-2}\to 0 \qquad \text{ as } k\to \infty.
\end{align*}
Thus \(\Gamma_1\) is not \(\varphi_\delta\)-separated.

On the other hand, let
\[\Gamma_2=\{a_k,\ a_k+y_k:\ k^{4\delta}\geq  2\},\] where $a_k= k-\delta$ and $y_k=k^{-4\delta}$.
Then $|2a_k\pi-2(a_k+y_k)\pi|=2y_k\pi=2 \pi k^{-4\delta}\to 0,$ as $k\to \infty$. So \(\Gamma_2\) is not \(\varphi_\pi\)-separated. However, \(\Gamma_2\) is \(\varphi_\delta\)-separated.

Since, we have $a_k<a_k+y_k<a_{k+1}<a_{k+1}+y_{k+1}$. We begin by establishing a lower bound for the $\varphi_\delta$ distance between $a_k$ and $a_k+y_k$. Therefore, consider
\begin{align*}
	|\varphi_\delta(a_k+y_k)-\varphi_\delta(a_k)|&=\varphi_\delta(a_k+y_k)-\varphi_\delta(a_k)=\int_{a_k}^{a_k+y_k}\varphi_\delta'(x)dx\\
	&\ge \int_{a_k}^{a_k+y_k} \frac{y_k}{(x-a_k)^2+y_k^2}dx=\int_{0}^{y_k}\frac{y_k}{t^2+y_k^2}dt= \frac{\pi}{4}.
	\end{align*}
Next, estimating the lower bound for the $\varphi_\delta$ distance between $a_{k+1}$ and $a_k+y_k$. With the help of  $a_k+y_k<a_{k+1}-y_{k+1}$, write
\begin{align*}
	|\varphi_\delta(a_{k+1})-\varphi_\delta(a_k+y_k)|&=\varphi_\delta(a_{k+1})-\varphi_\delta(a_k+y_k)=\int_{a_k+y_k}^{a_{k+1}}\varphi_\delta'(x)dx\\
	&\ge \int_{a_{k+1}-y_{k+1}}^{a_{k+1}}\varphi_\delta'(x)dx\ge \int_{a_{k+1}-y_{k+1}}^{a_{k+1}}\frac{y_{k+1}}{(x-a_{k+1})^2+y_{k+1}^2}dx\\
	&= \int_{-y_{k+1}}^0 \frac{y_{k+1}}{t^2+y_{k+1}^2}dt=\frac{\pi}{4}.
\end{align*}
Therefore, for any pair of points $\gamma_1,\gamma_2$ from the sequence $\Gamma_2$, one can show that $|\varphi_\delta(\gamma_1)-\varphi_\delta(\gamma_2)|\ge \frac{\pi}{4}$. Thus, \(\Gamma_2\) is \(\varphi_\delta\)-separated. Hence, the two
separation notions are not comparable.
\end{proof}
We now identify a general class of subsets of the real line on which the two notions of phase separation agree. The main idea is that each component interval becomes small in both phase metrics and that different intervals remain uniformly separated in both metrics.
\begin{proposition}\label{region}
Let \(\varphi_\pi(x)=2\pi x\), and let \(\varphi_\delta\) be the phase function associated with \(E_\delta\). Suppose that
\begin{equation}\label{region_B}
\mathcal B=\bigcup_{j\geq 1} J_j\subset\mathbb R
\end{equation}
is a union of compact intervals satisfying
$\operatorname{diam}_{\varphi_\pi}(J_j)\to 0$ and $\operatorname{diam}_{\varphi_\delta}(J_j)\to 0;$ where $\displaystyle\operatorname{diam}_{\varphi}(J_j) :=\sup_{x,y\in J_j} |\varphi(x)-\varphi(y)|$,
and \[\inf_{j\neq k}\operatorname{dist}_{\varphi_\pi}(J_j,J_k)>0,\qquad\inf_{j\neq k}\operatorname{dist}_{\varphi_\delta}(J_j,J_k)>0.
\]Then, for every sequence \(\Gamma\subset\mathcal B\), $\Gamma \text{ is }\varphi_\pi\text{-separated}\text{ if and only if }\Gamma\text{ is } \varphi_\delta\text{-separated}.$
\end{proposition}

\begin{proof}
	Let $\Gamma$ is $\varphi_\pi$-separated. This implies there exists a constant $\alpha>0$ such that 
	\begin{equation}\label{seperation}
		|\varphi_\pi(\gamma)-\varphi_\pi(\gamma')|\ge \alpha\qquad \text{ for all } \gamma,\gamma'\in\Gamma.
	\end{equation}
	Since $\operatorname{diam}_{\varphi_\pi}(J_j)\to 0$, therefore for every $\epsilon>0$, there exists some $n_0\in \mathbb N$ such that 
	\[\operatorname{diam}_{\varphi_\pi}(J_j)<\epsilon\qquad \forall ~~ j\ge n_0.\]
This holds for $\epsilon=\alpha$, and there exists some $j_0\in \mathbb N$ such that $\operatorname{diam}_{\varphi_\pi}(J_j)<\alpha$ for all $ j\ge j_0$. Assume that some $J_j$ with $j\ge j_0$ contains two distinct points $\gamma_0,\gamma_1\in \Gamma$, this gives
\[|\varphi_\pi(\gamma_0)-\varphi_\pi(\gamma_1)|\leq \operatorname{diam}_{\varphi_\pi}(J_j) <\alpha.\]
This contradicts equation \eqref{seperation}. Therefore, the number of elements of $\Gamma$ in any of $J_j$ such that $j\ge j_0$ is atmost one.

Now we consider $J_1,J_1,\ldots,J_{{j_0}-1}$. Since, $\Gamma$ is $\varphi_\pi$-seperated. Thus from equation \eqref{seperation} it follows that,
\[|\gamma-\gamma'|\ge \frac{\alpha}{2\pi}>0 \qquad \forall \gamma,\gamma'\in \Gamma.\]
Thus, $\Gamma $ is Euclidean seperated sequence and it can have only finitely many points in compact interval. Thus, $\displaystyle\Gamma\cap \left(\cup_{j=1}^{j_0-1}J_j\right)$ is finite. Because $\varphi_\delta$ is an increasing function, the $\varphi_\delta$-distance between distinct points of this finite set are positive. Hence, minimum of those distances, whenever there are atleast two points, is also positive. In particular, there exists $c_1>0$ such that
\begin{equation}\label{finiteintervals}
	|\varphi_\delta(\gamma)-\varphi_\delta(\gamma')|\ge c_1,
\end{equation}
whenever $\gamma,\gamma'\in \displaystyle\Gamma\cap \left(\cup_{j=1}^{j_0-1}J_j\right) .$	
	
	To show $\Gamma$ is $\varphi_\delta-$seperated. Take $\gamma_1,\gamma_2\in \Gamma.$ If they belong to different intervals then there exists some $c_2>0$ such that
	\[|\varphi_\delta(\gamma_1)-\varphi_\delta(\gamma_2)|\ge c_2.\]
	If $\gamma_1$ and $\gamma_2$ belong to same interval then that interval must be among $J_1,J_2,\ldots,J_{{j_0}-1} $ and it follows from equation \eqref{finiteintervals} that \[|\varphi_\delta(\gamma_1)-\varphi_\delta(\gamma_2)|\ge c_1.\]Take $\min\{c_1,c_2\}=c'$, this yields $|\varphi_\delta(\gamma_1)-\varphi_\delta(\gamma_2)|\ge c'.$ Consequently, $\Gamma$ is $\varphi_\delta$-seperated.  
	
	Conversely assume $\Gamma$ is $\phi_\delta$-seperated. Therefore, there exists a positive constant $\beta>0$ such that 
	\begin{equation}\label{seperationconstant2}
		|\varphi_\delta(\gamma)-\varphi_{\delta}(\gamma')|\ge \beta \qquad \text{ for all }\gamma,\gamma'\in\Gamma.
	\end{equation}
	Similar to the previous case there exists some $j_1$ such that the number of elements of $\Gamma$ in any of $J_j$ such that $j\ge j_1$ is atmost one. Since \(\varphi_\delta\) is continuous and strictly increasing, the image of each compact interval \(J_i; ~1\leq i \leq j_1-1\), under \(\varphi_\delta\) is a bounded interval of finite length. 
	Hence, with the help of equation \eqref{seperationconstant2} one can say that each of the finitely many intervals \(J_1,\ldots,J_{j_1-1}\) contains only finitely many points of \(\Gamma\). Hence their union contains only finitely many points of \(\Gamma\). Since \(\varphi_\pi(x)=2\pi x\) is strictly increasing, the same argument as above leads to $\Gamma$ is	$\varphi_\pi$-seperated.
\end{proof}
The following example provides a region in which the two notions of separation are equivalent. A natural question is whether this region contains sequences whose density satisfies the hypotheses of the main result in \cite{MarzoNitzanOlsen2012}. The answer is affirmative, such sequences do indeed exist within this region.

\begin{example}
For $|n|\ge N$, define
	\[\lambda_n=n-\operatorname{sgn}(n)\delta ,\qquad\varepsilon_n=|n|^{-8\delta}. \]
	where  \(N^{8\delta}>4\), and
	\[P_n=[\lambda_n,\lambda_n+\varepsilon_n], \qquad Q_n=\left[\lambda_n+\frac12, \lambda_n+\frac12+\varepsilon_n\right],\]
	and set
	\[\mathcal B_N=\bigcup_{|n|\geq N}(P_n\cup Q_n).\]
	Then, for every sequence \(\Gamma\subset\mathcal B_N\),
	$\Gamma \text{ is } \varphi_\delta\text{-separated} $
	 if and only if $\Gamma \text{ is }\varphi_\pi\text{-separated}.$
	Moreover, \(\mathcal B_N\) contains separated sequences whose lower and upper \(\varphi_\delta\)-densities lie on either side of \(1/2\pi\).
	
\end{example}
\begin{proof}
	The phase derivative corresponding to phase function $\varphi_\delta$ is given by
	\[ \varphi_\delta'(x) = \frac{2}{x^2+1}+\sum_{j\geq 1}\left[\frac {2j^{-4\delta}}{(x- a_j)^2+j^{-8\delta}}+\frac{2j^{-4\delta}}{(x+a_j)^2+j^{-8\delta}}\right],\]
	where $a_j=j-\delta$. The function \(\varphi_\delta\) is strictly increasing on \(\mathbb R\).
	
	We first prove that the intervals \(P_n\) and \(Q_n\) have small diameter in
	both metrics. Let $n_0>N$ be arbitrary, 
	\[\operatorname{diam}_{\varphi_\delta}(P_{n_0}) =\sup_{x,y\in P_{n_0}} |\varphi_\delta(x)-\varphi_\delta(y)|,\]
	Since $\varphi_\delta$ is a strictly increasing function, \[\operatorname{diam}_{\varphi_\delta}(P_{n_0})= \varphi_\delta(\lambda_{n_0} + \varepsilon_{n_0} ) -\varphi_\delta(\lambda_{n_0})=\int_{\lambda_{n_0}}^{\lambda_{n_0} + \varepsilon_{n_0}}\varphi_\delta'(x)dx.\]
	Since, for $n_0>N$, $a_{j_0}$ and $\lambda_{j_0}$ will be same. So, we divide the above integral into three parts as 
	\[\int_{\lambda_{n_0}}^{\lambda_{n_0}+\varepsilon_{n_0}}\frac{1}{x^2+1}dx,\quad \sum_{\substack{j\in\mathbb Z\setminus\{0\}\\j\neq n_0}}\int_{\lambda_{n_0}}^{\lambda_{n_0} + \varepsilon_{n_0}}\frac {|j|^{-4\delta}}{(x- \lambda_j)^2+|j|^{-8\delta}}dx ,\quad\text{and } \int_{\lambda_{n_0}}^{\lambda_{n_0} + \varepsilon_{n_0}}\frac {{n_0}^{-4\delta}}{(x- \lambda_{n_0})^2+{n_0}^{-8\delta}}dx.\]
	The contribution of the \(n_{0}\)-th Poisson kernel to the
	\(\varphi_\delta\)-length of \(P_{n_0}\) is
	\begin{equation}\label{for n_0}
		\int_{\lambda_{n_0}}^{\lambda_{n_0}+\varepsilon_{n_0}}\frac{{n_0}^{-4\delta}}{(x-\lambda_{n_0})^2+{n_0}^{-8\delta}}dx =	\int_0^{\varepsilon_{n_0}}  \frac{{n_0}^{-4\delta}}{t^2+{n_0}^{-8\delta}} dt=\arctan\left(\frac{{n_0}^{-8\delta}}{{n_0}^{-4\delta}}\right)=\arctan({n_0}^{-4\delta}),
	\end{equation}
and also, 
\begin{equation}\label{for1/1+x^2}
	\int_{\lambda_{n_0}}^{\lambda_{n_0} + \varepsilon_{n_0}}\frac{1}{x^2+1}dx \leq \varepsilon_{n_0}= |n_0|^{-8\delta}.
	\end{equation}
Now we estimate the upper bound for the summation term over $j\neq n_{0}$. Since, $|{n_{0}}|^{8\delta}>N^{8\delta}>4$, we have $|x-\lambda_{{n_{0}}}|<\frac14$, this gives \[|x-\lambda_j|\ge |\lambda_j-\lambda_{n_{0}}|-|x-\lambda_{n_{0}}|\ge |\lambda_j-\lambda_{n_{0}}|-\frac14 \ge \frac34|\lambda_j-\lambda_{n_{0}}| , \quad j\neq n_0.\]
Using this we get,
\begin{equation}\label{fornotequaltoj}
	\sum_{\substack{j\in\mathbb Z\setminus\{0\}\\j\neq n_0}}\frac {|j|^{-4\delta}}{(x- \lambda_j)^2+|j|^{-8\delta}}\leq\sum_{\substack{j\in\mathbb Z\setminus\{0\}\\j\neq n_0}}\frac {1}{(x- \lambda_j)^2}\leq \sum_{\substack{j\in\mathbb Z\setminus\{0\}\\j\neq n_0}}\frac {16}{9|\lambda_j-\lambda_{n_{0}}|^2}\leq 2 \sum_{j=1}^{\infty}\frac{1}{j^2}<C . 
\end{equation}
By equations \eqref{for n_0}, \eqref{for1/1+x^2} and \eqref{fornotequaltoj} we get,
\[\int_{\lambda_{n_0}}^{\lambda_{n_0} + \varepsilon_{n_0}}\varphi_\delta'(x)dx\leq 2|n_0|^{-8\delta}+2C |n_0|^{-8\delta}+2\arctan({n_0}^{-4\delta}) . \]
Since $n_0$ is arbitrary  $n_0>N$. Therefore, one can write above expression for any $n>N$ and as $n$ tends to infinity we get, 	
\[\int_{\lambda_{n}}^{\lambda_{n} + \varepsilon_{n}}\varphi_\delta'(x)dx\leq 2|n|^{-8\delta}+2C |n|^{-8\delta}+2\arctan({n}^{-4\delta}) \to 0 \quad \text{ as } n\to\infty. \]
Therefore,
\[\operatorname{diam}_{\varphi_\delta}(P_n)= \varphi_\delta(\lambda_n+\varepsilon_n)-\varphi_\delta(\lambda_n)\longrightarrow 0.\]
	Also,
	\[\operatorname{diam}_{\varphi_\pi}(P_n)=2\pi \varepsilon_n\longrightarrow 0.\]
Similarly, \(Q_n\) stays a distance comparable to \(1/2\) away from the
	nearest points \(\lambda_n\) and \(\lambda_{n+1}\). Indeed, for large $n$ when $x\in Q_n$, 
	\[|x-\lambda_n|\ge \frac12, \qquad |x-\lambda_{n+1}|=|x-(\lambda_n+1)|\ge \frac14.\]
	 Notice that between any interval of length one such that $j+\frac13\leq |x-\lambda_j|\leq \frac43, ~~ x\in Q_n,$ there can be atmost two $\lambda_j$ because $\lambda_j$'s are separated by distance 1. Therefore, on dividing the summation over $j\neq n_{0}$ into such intervals we obtain the following estimates   
	\[\sum_{j\in\mathbb Z\setminus\{0\}}\frac {|j|^{-4\delta}}{(x- \lambda_j)^2+|j|^{-8\delta}}dx \leq \sum_{j\neq0}\frac {1}{(x- \lambda_j)^2}dx \leq 2 \sum_{j=1}^{\infty} \frac{1}{(j+1/3)^2}\leq C,\]
	since $\sum_{j=1}^{\infty} \frac{1}{(j+1/3)^2}$ is a convergent series.
	Additionally, $\frac{1}{1+x^2}\leq 1$. Hence
	\(\varphi_\delta'\) is uniformly bounded on \(Q_n\), and so
	\[\operatorname{diam}_{\varphi_\delta}(Q_n) \leq C\varepsilon_n\longrightarrow 0. \]
	One can also show that
	\[ \operatorname{diam}_{\varphi_\pi}(Q_n) = 2\pi\varepsilon_n
	\longrightarrow 0. \]
	
Next we prove that distinct intervals are separated by a fixed positive
	distance in both metrics. Clearly,
	\[\operatorname{dist}(P_n,Q_n)\geq \frac14, \qquad \operatorname{dist}(Q_n,P_{n+1})\geq \frac14. \]
	Thus distinct intervals are uniformly separated in the Euclidean metric,
	and therefore also in the \(\varphi_\pi\)-metric, since
	$ \varphi_\pi(x)=2\pi x. $\\
	We now check the \(\varphi_\delta\)-metric. The smallest possible phase distance occurs between the right
	endpoint of \(P_n\) and the left endpoint of \(Q_n\). Hence
	\begin{align*}
	\operatorname{dist}_{\varphi_\delta}(P_n,Q_n)=\int_{\lambda_n+\varepsilon_n}^{\lambda_n+1/2} \varphi_\delta'(x)dx &\geq
	\int_{\lambda_n+\varepsilon_n}^{\lambda_n+1/2}
	\frac{{n}^{-4\delta}}{(x-\lambda_n)^2+{n}^{-8\delta}}dx= \int_{\varepsilon_n}^{1/2}\frac{{n}^{-4\delta}}{t^2+{n}^{-8\delta}}dt\\
	&=\arctan\left(\frac{1}{2{n}^{-4\delta}}\right)- \arctan\left(\frac{\varepsilon_n}{{n}^{-4\delta}}\right)\\
	&=\arctan\left(\frac{1}{2{n}^{-4\delta}}\right)-\arctan({n}^{-4\delta})>\frac{\pi}{4}-\arctan(1/2).
	\end{align*}
 Thus, there exists a constant \(c_0>0\) such that $\operatorname{dist}_{\varphi_\delta}(P_n,Q_n)\geq c_0$,
	for all \(|n|\geq N\). Similarly, using the Poisson kernel centered at \(\lambda_{n+1}\), one
	obtains $\operatorname{dist}_{\varphi_\delta}(Q_n,P_{n+1})\geq c_1>0$. Hence all distinct intervals in the family $\{P_n,Q_n: |n|\geq N\}$ are separated by a fixed positive distance in the \(\varphi_\delta\)-metric.
	
Since all conditions of Proposition \ref{region} are satisfied here for $P_n$ and $Q_n$, every seperated sequence \(\Gamma\subset\mathcal B_N\), $\Gamma \text{ is }\varphi_\pi\text{-separated}\text{ if and only if }\Gamma\text{ is } \varphi_\delta\text{-separated}.$
	
	It remains to prove the density claims. Theorem \ref{densities} establishes that both upper and lower densities corresponding to phase function $\varphi_\delta$ are same as the densities associated with $\varphi_\pi$. Additionally, the densities corresponding to 
	$\varphi_\pi$ are comparable with Euclidean densities because $\varphi_\pi$ is nothing but $2\pi x$. Therefore, it is easy to verify that sequence
\[\Gamma_1=\left\{ \lambda_n,\lambda_n+\frac12:\ |n|\geq N	\right\}\subset\mathcal B_N,\]
	has upper and lower density corresponding to $\varphi_\delta$ to be $1/\pi>1/2\pi$. On the other hand, on taking
	\[\Gamma_2=\{\lambda_{2n}: |2n|\geq N\}\subset\mathcal B_N,\]
	we get,
	\[D_{\varphi_\delta}^-(\Gamma_2)= D_{\varphi_\delta}^+(\Gamma_2) = \frac{1}{4\pi} < \frac{1}{2\pi}. \]
	Thus inside \(\mathcal B_N\), one can find separated sequences whose lower
	and upper \(\varphi_\delta\)-densities are greater than \(1/2\pi\), and also
	separated sequences whose lower and upper \(\varphi_\delta\)-densities are
	less than \(1/2\pi\). This completes the proof.
\end{proof}
We are now in a position to state the analogue of the main result from \cite{MarzoNitzanOlsen2012}. The notation \(\mathcal{B}\) used below is defined in \eqref{region_B}.

\begin{theorem}\label{sampling}
Let $E_\delta $ be a Hermite-Biehler function with phase function $\varphi_{\delta}$ which is a non-doubling measure as defined in Lemma \ref{non_doubling_measure}. Every real sampling sequence for $H(E_\delta)$, $\Gamma\subset \mathcal B$, is a finite union of $\varphi_\delta$-separated sequences. Moreover, there exists a $\varphi_\delta$-separated $\Gamma'\subset \Gamma$ which is sampling and satisfies $D^-_{\varphi_\delta}(\Gamma')\ge 1/2\pi$. Conversely, every real sequence $\Gamma\subset \mathcal B$ which is a finite union of $\varphi_\delta$-separated sequences $\Gamma'\subset\Gamma$ such that $D^-_{\varphi_\delta}(\Gamma')> 1/2\pi$ is sampling.
\end{theorem}
\begin{proof}
	The proof of the theorem is obtained by combining Theorem \ref{densities} and Proposition \ref{region} with the equation \eqref{equivalences} in Subsection \ref{equivalence}. Indeed, using equation \eqref{equivalences}, every real sampling sequence for $H(E_\delta)$ is also a real sampling sequence for $H(e^{-i\pi z})$ and using Theorem 1.4 of \cite{MarzoNitzanOlsen2012}, every real sampling sequence $\Gamma \subset \mathcal B$ is a finite union of $\varphi_\pi$-seperated sequences. Since, $\varphi_\pi$-seperation  and $\varphi_\delta$-seperation are equivalent due to  Proposition \ref{region}, this gives every real sampling sequence for $H(E_\delta)$, $\Gamma\subset \mathcal B$, is a finite union of $\varphi_\delta$-separated sequences. Similarly, the rest of the proof follows by similar argument.
	\end{proof}
\begin{theorem}\label{interpolating}
	Let $E_\delta $ be a Hermite-Biehler function with phase function $\varphi_{\delta}$ which is a non-doubling measure as defined in Lemma \ref{non_doubling_measure}. Every real interpolating sequence  $\Gamma\subset \mathcal B$ for $H(E_\delta)$ is $\varphi_\delta$-separated with $D^+_{\varphi_\delta}(\Gamma)\leq 1/2\pi$. Conversely, every real sequence $\Gamma\subset \mathcal B$ which is $\varphi_\delta$-separated with $D^+_{\varphi_\delta}(\Gamma)< 1/2\pi$ is interpolating.
	
\end{theorem}
\begin{proof}
	The proof of the theorem is obtained by combining Theorem \ref{densities} and Proposition \ref{region} with the equation \eqref{equivalences} in Subsection \ref{equivalence}.
\end{proof}
We conclude this article with a few open problems arising from the results presented in Section \ref{sec4}.\\
	\textbf{Problem 1.} Extend the class of functions represented by $E_\delta$, as defined in Lemma \ref{non_doubling_measure}, and characterize a broader class of functions for which Theorems \ref{sampling} and \ref{interpolating} remain valid.\\
	\textbf{Problem 2.} Determine whether the class of regions considered in Theorems \ref{sampling} and \ref{interpolating} can be enlarged while preserving the conclusions of these theorems.
\section*{Declarations}

\subsection*{Conflicts of Interest}

The authors declare that there are no known competing financial or personal interests that could have appeared to influence the work reported in this paper.


\begin{thebibliography}{9}
	

	\bibitem{Baranov2000}
	A.~D. Baranov,
	\textit{Differentiation in the Branges spaces and embedding theorems},
	J. Math. Sci. (New York) \textbf{101} (2000), no.~2, 2881--2913.
	Nonlinear equations and mathematical analysis.
	
	\bibitem{deBranges1968}
	L. de Branges,
	\textit{Hilbert Spaces of Entire Functions},
	Prentice-Hall, Englewood Cliffs, NJ, 1968.
	
	\bibitem{fuhr2007sampling}
	H. Fuhr and K. Gr\"ochenig,
	\textit{Sampling theorems on locally compact groups from oscillation estimates,}
	Mathematische Zeitschrift, 255 (2007), 177--194.
	
	\bibitem{fuhr2017density}
	H. Fuhr, K. Gr\"ochenig, A. Haimi, A. Klotz, and J. L. Romero,
	\textit{Density of sampling and interpolation in reproducing kernel Hilbert spaces,}
	Journal of the London Mathematical Society, 96 (2017), no. 3, 663--686.
	
	\bibitem{mashregi}
	S. Garcia, J. Mashreghi, and W. T. Ross.
	{\em Introduction to model spaces and their operators}, Vol. 148.
	Cambridge University Press (2016).
	
	\bibitem{MarzoNitzanOlsen2012}
	J.~Marzo, S.~Nitzan, and K.~Seip Olsen,
	\textit{Sampling and interpolation in de Branges spaces with doubling phase},
	J. Analyse Math. \textbf{117} (2012), 365--395.
	
	\bibitem{marakov2005}
	N.~Makarov and A.~Poltoratski,
	\textit{Meromorphic inner functions, Toeplitz kernels and the uncertainty principle},
	in \textit{Perspectives in Analysis: Essays in Honor of Lennart Carleson's 75th Birthday},
	Springer, Berlin, Heidelberg, 2005.
	
	\bibitem{Mashreghi_representaion_theory}
	J. Mashreghi,
	\textit{Representation Theorems in Hardy Spaces},
	Cambridge University Press, Cambridge, 2009.
	
	\bibitem{Mitkovski2010}
	M. Mitkovski, and A. Poltoratski.
	{\em P\'olya sequences, Toeplitz kernels and gap theorems.}
	Advances in Mathematics 224.3 (2010) 1057-1070.
	
	\bibitem{Nikolskii1986}
	N.~K. Nikolskii.
	{\em Treatise on the Shift Operator},
	Springer-Verlag, Berlin--New York, 1986.
	
	\bibitem{Nikolskii2002I}
	N.~K. Nikolskii,
	\textit{Operators, Functions, and Systems: An Easy Reading, Vol.~I},
	American Mathematical Society, Providence, RI, 2002.
	
	\bibitem{Nikolskii2002II}
	N.~K. Nikolskii,
	\textit{Operators, Functions, and Systems: An Easy Reading, Vol.~II},
	American Mathematical Society, Providence, RI, 2002.
	
	\bibitem{Levin1996}
	B.~Ya. Levin,
	\textit{Lectures on Entire Functions},
	Translations of Mathematical Monographs, vol.~150,
	American Mathematical Society, Providence, RI, 1996.
	In collaboration with and with a preface by Yu.~Lyubarskii,
	M.~Sodin, and V.~Tkachenko; translated from the Russian manuscript
	by V.~Tkachenko.
	
	\bibitem{LyubarskiiSeip2002}
	Yurii I. Lyubarskii and Kristian Seip,
	\textit{Weighted Paley--Wiener spaces},
	J. Amer. Math. Soc. \textbf{15} (2002), no.~4, 979--1006.
	
	\bibitem{landau1967necessary}
	H. J. Landau,
	\textit{Necessary density conditions for sampling and interpolation of certain entire functions,}
	Acta Mathematica, 1967.
	
	\bibitem{seip1992densityFockI}
	K. Seip,
	\textit{Density theorems for sampling and interpolation in the Bargmann--Fock space I,}
	Journal f\"ur die reine und angewandte Mathematik, 429 (1992), 91--106.
	
	\bibitem{Seip2004}
	K. Seip,
	\textit{Interpolation and Sampling in Spaces of Analytic Functions},
	American Mathematical Society, Providence, RI, 2004.
	
	
	
	
	
\end{thebibliography}
	\end{document}